# Corrective Control: Stability Analysis of Unified Controller Combining Frequency Control and Congestion Management

O.O. Khamisov*, T. S. Chernova*, J.W. Bialek*, S. H. Low $^{\&}$
* Center for Energy Systems, Skolkovo Institute of Science and Technology (Skoltech), Moscow, Russia
$^{\&}$ California Institute of Technology, Pasadena, USA

## Summary

This paper analyses stability of the Unified Controller (UC) that combines frequency control and congestion management and therefore makes it possible to move from preventive to corrective power system control. Earlier work by the authors of UC proved asymptotic stability of the methodology but the proof was based on a simplified first-order model of the turbine and turbine governor. We show that a higher order model of the turbine governor leads to eigenvalues with small but positive real parts. Consequently, we develop a modification of the methodology that decouples the physical and control systems and therefore results in all the eigenvalues having negative real parts. We illustrate the effectiveness of the modification on a realistic model of 39-bus model of New England power system implemented in Power System Toolbox (PST).

## 1  Notation and nomenclature

Let $\mathbb{R}$ be set of real numbers and $\mathbb{N}$ the set of natural numbers. For finite set $X$ we denote its cardinality as $|X|$. For arbitrary vector $x$ of size $n$, and set $I \subseteq \{1, ..., n\}$ we use $x_I = (x_i, i \in I)$ to denote its subvector.

*Sets*
| | |
|---|---|
| $G$ | Set of generator buses. |
| $L$ | Set of load buses. |
| $N = G \cup L, |N| = n$ | Set of all buses. |
| $\mathcal{E} \subseteq N \times N, |\mathcal{E}| = m$ | Set of lines. |

*Variables*

$\theta_i, i \in N$      Rotor angle (for generator buses) or bus voltage angle (for load buses) with respect to a synchronously rotating reference axis

$\omega_i, i \in N$      Rotor speed deviation (for generator buses) or bus frequency deviation (for load buses) from the nominal value

$P_{ij}, ij \in \mathcal{E}$      Lossless line power, $P_{ij} = -P_{ji}$.

$p_i, i \in N$      Power control command at generator and load nodes

$p_i^M, i \in G$      Mechanical power injection at generators.

$v_i, i \in G$      Turbine valve position

*Parameters*

$M_i, i \in G$      Inertia constants of the generators.

$D_i \omega_i, i \in N$      The aggregate frequency sensitivity of generators and loads.

$r_i, i \in N$      Unknown disturbance modeled as simultaneous step changes in power injections on an arbitrary subset of the buses

$B_{ij}, ij \in \mathcal{E}$      Line susceptance

$T_i^t, i \in G$      Time constant of turbine

$T_i^g, i \in G$      Time constant of governor.

## 2  Introduction

The essence of power system control is to maintain system security at minimal cost. This paper deals with arguably two most important aspects of system security: keeping frequency within tight bounds around the nominal value, which is traditionally executed by centralized Automatic Generation Control (AGC), and maintaining power flows below line limits so that the lines are not overloaded (congestion management). The latter is enforced by finding optimal (i.e. minimum cost) dispatch using Security Constrained Optimal Power Flow (SCOPF). Typically (N-1) security standard is observed which means that no single contingency (i.e. a loss of a single element such as a line or a generator) should result in such redistribution of power flows that any line is overloaded. That philosophy of system control, combining real-time AGC control with SCOPF run at certain intervals (usually between 5 mins and 1 hour), has evolved over decades and is considered almost sacrosanct by the industry. It has proved to be effective as large system blackouts are very rare. However the excellent system security comes at a considerable extra cost. The reason is that SCOPF is a tool for preventive control, i.e. finding an operating point that is safe *should* a disturbance happen and therefore requires carrying a significant reserve at all times. This problem has been recognized over the last decade or so and research has been undertaken to move from preventive to corrective control, i.e. to a control that maintains system security *when* a disturbance happens. Corrective control typically reduces the requirement to keep a costly reserve because if any contingency happens, a remedial control is activated to return the system to a secure operation.

In this paper we utilize the concept of Unified Control (UC) [1,2] that combines frequency control with congestion management and therefore makes it possible to relieve transmission constraints when a disturbance (i.e. a power balance change or a line trip) happens. Every control agent participating in UC (a generator, controllable load or energy storage) measures its local frequency and power flows, performs moderate computations, and communicates with its neighbouring agents in a communication network with an arbitrary topology (as long as it connects all the agents).

UC, when compared to standard Automatic Generation Control (AGC), has many advantages [1,2]. It is decentralized, it can be used not only for generators but also for controllable loads and therefore be used in future low-inertia systems, and it can be run in an open-architecture in a plug-and-play manner - i.e. any additional agent (generator or controllable load) can connect to the system and participate in frequency control without coordinating with System Operator. AGC requires communication of each generator with System Operator while UC requires each agent to communicate with its neighbors and exchange the values of a few variables. This makes it possible for a large number of additional distributed resources to participate in frequency control and relieving transmission constraints without overwhelming the System Operator. Cost functions of individual agents are private to agents and do not need to be shared. The ability of UC to control frequency while keeping power flows within their limits offers an excellent opportunity to implement corrective, rather than preventive, control. With UC, the optimal dispatch could be then calculated using straight OPF rather than SCOPF (i.e. without considering contingencies or with a reduced set of contingencies) as any contingencies would be dealt with by UC in real-time when they happen.

The main aim of this paper is to analyze stability of UC. While the earlier work by the authors of the algorithm [1,2] proved asymptotic stability of the methodology, the proof was based on a simplified first-order model of the turbine and turbine governor. We show that a higher-order model of the turbine governor leads to eigenvalues with small but positive real parts leading to sustained oscillations. Consequently we develop a modification of the methodology that results in all the eigenvalues having negative real parts in the higher-order model. Testing the modification on IEEE 39-bus model containing 10 generators and using realistic dynamic models of the generators, exciters and turbine governors implemented in PST has proved effectiveness of the methodology.

## 3 Power network model

We consider a generator dynamic model combined with speed governor and turbine model [4]. We assume that all voltages are constant and equal to 1 p.u, line resistances are much smaller than the reactances, and reactive powers are neglected. Dynamics of the assumed power system model is defined by the system of differential algebraic equations.

$$\dot{\theta}_i = \omega_i, i \in N \quad (1a)$$

$$M_i \dot{\omega}_i = -D_i \omega_i - \sum_{ij \in \mathcal{E}} P_{ij} + p_i^M + r_i, i \in G, \quad (1b)$$

$$0 = -D_i \omega_i - \sum_{ij \in \mathcal{E}} P_{ij} + p_i + r_i, i \in L, \quad (1c)$$

$$P_{ij} = B_{ij} \sin(\theta_i - \theta_j), ij \in \mathcal{E}, \quad (1d)$$

$$T^t \dot{p}_i^M = -p_i^M + v_i, i \in G, \quad (1e)$$

$$T^g \dot{v}_i = -v_i + p_i, i \in G. \quad (1f)$$

Here equations (1b) express generator dynamics, (1c) are power balance equations for load buses, (1d) are equations defining lossless nonlinear power flows, equations (1e) and (1f) describe turbine and governor dynamics, respectively, as 1st-order lags.

Our aim is to derive a control $p_i, i \in N$, that delivers system (1) into an equilibrium $(\theta_i^*, \omega_i^*, P_{ij}^*, (p_i^M)^*, v_i^*)$ that satisfies a number of constraints.

Power flows must be kept within the line limits (congestion management):

$$\underline{P}_{ij} \leq P_{ij}^* \leq \overline{P}_{ij}, ij \in \mathcal{E}. \quad (2)$$

where $\underline{P}_{ij}, \overline{P}_{ij}$ are the lower and upper line limits, respectively.

Inter-area line flows must be restored to scheduled values:

$$\sum_{ij \in A_k} P_{ij}^* = P_k^{area}, \quad (3)$$

Here $A_k \subseteq \mathcal{E}, k \in \Gamma$, are sets of lines that connect different areas, i.e., $i$ is in area $k$, or $j$ is in area $k$, but not both. $P_k^{area}, k \in \Gamma$, are scheduled net flows between areas $k$ and the rest of the network.

Control must be within control limits at any moment of time:

$$\underline{p}_i \leq p_i(t) \leq \overline{p}_i, i \in N, t \geq 0. \quad (4)$$

The aim of frequency control is to restore the nominal frequency at the equilibrium:

$$\omega_i^* = 0, i \in N. \quad (5)$$

## 4 Unified Control (UC)

Full description and derivation of the Unified Control can be found in [1, 2] so we only provide a short review here. The aim of UC is to deliver the system to an equilibrium that minimizes the deviation from the pre-disturbance dispatch (which is assumed to be economically optimal) while satisfying constraints (2)-(4). The UC controller is described by the set of differential equations:

$$\hat{P}_{ij} = B_{ij} \sin(\phi_i - \phi_j), ij \in \mathcal{E}, \quad (6a)$$

$$\dot{\lambda}_i = K_i^\lambda \left( M_i \dot{\omega} + D_i \omega_i + \sum_{ij \in \mathcal{E}} (P_{ij} - \hat{P}_{ij}) \right), \quad (6b)$$

$$i \in G,$$

$$\dot{\lambda}_i = K_i^\lambda \left( D_i \omega_i + \sum_{ij \in \mathcal{E}} (P_{ij} - \hat{P}_{ij}) \right), \quad (6c)$$

$$i \in L,$$

$$\dot{\phi}_i = K_i^\phi \sum_{ij\in\mathcal{E}} B_{ij}(\lambda_i - \lambda_j - \rho_{ij}^+ \qquad (6d)$$
$$+\rho_{ij}^- - \sum_{k:ij\in A_k} \pi_k \bigg), i \in N$$
$$\dot{\rho}_{ij}^+ = K_{ij}^{\rho^+}[\hat{P}_{ij} - \overline{P}_{ij}]_{\rho_{ij}^+}^+, ij \in \mathcal{E}, \qquad (6e)$$
$$\dot{\rho}_{ij}^- = K_{ij}^{\rho^-}[\underline{P}_{ij} - \hat{P}_{ij}]_{\rho_{ij}^-}^+, ij \in \mathcal{E}, \qquad (6f)$$
$$\dot{\pi}_k = K_k^\pi \bigg(\sum_{ij\in A}\hat{P}_{ij} - P_k^{area}\bigg), k \in \Gamma, \qquad (6g)$$
$$p_i = -\alpha_i(\omega_i + \lambda_i), i \in N. \qquad (6h)$$

Here $K_i^\lambda, K_i^\phi, K_{ij}^{\rho^+}, K_{ij}^{\rho^+}$ are control gains. Operator $[x]_y^+$ for any $x, y \in \mathbb{R}$ is defined as follows:

$$[x]_y^+ = \begin{cases} x, & \text{if } y > 0 \text{ or } x > 0, \\ 0, & \text{otherwise.} \end{cases} \qquad (7)$$

Note that (6h) effectively derives the control as the sum of two control signals, one proportional to $\omega_i$ and the other to $\lambda_i$, respectively. The first one is equivalent to the droop control so by making $\alpha_i = 1/R_i$, where $R_i$ is the effective gain of the governing system, we can implement standard droop control. The second signal is equivalent to the secondary frequency control that removes the frequency deviation.

The design rationale for the UC controller (6) is that the closed-loop system that consists of (a simplified version of) the physical system (1) and the cyber system (6) carries out a primal-dual algorithm in real time over the power network to solve the following optimization problem (see [1, 2] for details):

$$\min_{\tilde{\theta}_i, \tilde{P}_{ij}, \tilde{p}_i} \sum_{i \in N} \frac{1}{2}\alpha_i \tilde{p}_i^2 \qquad (8a)$$

subject to
$$r_i + \tilde{p}_i - \sum_{ij\in\mathcal{E}}\tilde{P}_{ij} = 0, i \in N, \qquad (8b)$$
$$\tilde{P}_{ij} = B_{ij}\sin(\tilde{\theta}_i - \tilde{\theta}_j), ij \in \mathcal{E}, \qquad (8c)$$
$$\underline{P}_{ij} \leq \tilde{P}_{ij} \leq \overline{P}_{ij}, ij \in \mathcal{E}, \qquad (8d)$$
$$\sum_{ij\in A_k}\tilde{P}_{ij} = P_k^{area}, \qquad (8e)$$
$$\underline{p}_i \leq \tilde{p}_i \leq \overline{p}_i, i \in N. \qquad (8f)$$

Here symbol $\sim$ is used to distinguish solutions of (1) and variables in (8). This optimization problem encodes the control goal of UC. Here $\alpha_i$ are disutility coefficients. Let us consider equilibrium of (1) with $(\theta_i^*, \omega_i^*, P_{ij}^*, (p_i^M)^*, v_i^*)$ with $\theta_i^* = \tilde{\theta}_i$, $P_{ij}^* = \tilde{P}_{ij}$ and $p_i = \tilde{p}_i$. From (1e), (1f) we get $(p_i^M)^* = v_i^* = \tilde{p}_i$. Substituting it to (1b), (1c) and using (8b) we get $\omega_i^* = 0$. Therefore equilibrium point given by the solution of optimization problem (8) satisfies (5). Constraints (2), (3) and (4) are explicitly present in (8) and therefore are also satisfied. Thus any control scheme that converges to the solution of (8) gives the desired control values. The Unified Control is derived as a variant of primal-dual algorithm to solve a *suitably modified version* of problem (8) – see [1-2] for details - and its Lagrangian dual. Therefore the generator dynamics automatically carry out the primal-dual algorithm over the network in real time.

## 5 Decoupled UC (DUC)

In [1], [2], stability of Unified Control (6) was proved for a system model that has only one turbine equation
$$T^t\dot{p}_i^M = -p_i^M + p_i, i \in G, \qquad (9)$$
and with linearized power flow equations. In practice, the turbine and its governing system is of a higher order [4]. In this paper we model it by two differential equations (1e), (1f) representing the turbine and turbine governor. We show below that, while the closed-loop (linearized) system with only one turbine equation [1]-[2] has stable eigenvalues, the closed-loop system under UC (6) with both (1e) and (1f) has eigenvalues with positive real parts and therefore is unstable. This motivates the following modification that attempts to decouple the control system and the physical system:

$$\hat{P}_{ij} = B_{ij}\sin(\phi_i - \phi_j), ij \in \mathcal{E}, \qquad (10a)$$
$$\dot{\lambda}_i = K_i^\lambda \bigg(M_i\dot{\omega}_i + D_i\omega_i + \sum_{ij\in\mathcal{E}}(P_{ij} - \hat{P}_{ij}) - \alpha_i\lambda_i \quad (10b)$$
$$- \tilde{p}_i^M\bigg), i \in G,$$
$$\dot{\lambda}_i = K_i^\lambda \bigg(D_i\omega_i + \sum_{ij\in\mathcal{E}}(P_{ij} - \hat{P}_{ij}) - \alpha_i\lambda_i \quad (10c)$$
$$- p_i\bigg), i \in L,$$
$$\dot{\phi}_i = K_i^\phi \sum_{ij\in\mathcal{E}} B_{ij}(\lambda_i - \lambda_j - \rho_{ij}^+ \qquad (10d)$$
$$+\rho_{ij}^- - \sum_{k:ij\in A_k}\pi_k\bigg), i \in N$$
$$\dot{\rho}_{ij}^+ = K_{ij}^{\rho^+}[\hat{P}_{ij} - \overline{P}_{ij}]_{\rho_{ij}^+}^+, ij \in \mathcal{E}, \qquad (10e)$$
$$\dot{\rho}_{ij}^- = K_{ij}^{\rho^-}[\underline{P}_{ij} - \hat{P}_{ij}]_{\rho_{ij}^-}^+, ij \in \mathcal{E}, \qquad (10f)$$
$$\dot{\pi}_k = K_k^\pi \bigg(\sum_{ij\in A}\hat{P}_{ij} - P_k^{area}\bigg), k \in \Gamma, \qquad (10g)$$
$$p_i = -\alpha_i(\omega_i + \lambda_i), i \in N, \qquad (10h)$$
$$\tilde{T}^t\dot{\tilde{p}}_i^M = -\tilde{p}_i^M + \tilde{v}_i, i \in G, \qquad (10i)$$
$$\tilde{T}^g\dot{\tilde{v}}_i = -\tilde{v}_i + p_i, i \in G. \qquad (10j)$$

**Design rationale.** Fig. 1 shows eigenvalues for the modified IEEE 39-bus New England system [5] using a linearized system model (1) but without equation (1e)

- this is similar to the model in [1-2]. In modeling, we have considered all nodes as generator nodes, assumed no congestion management and increased control gains to emphasize the points made below.

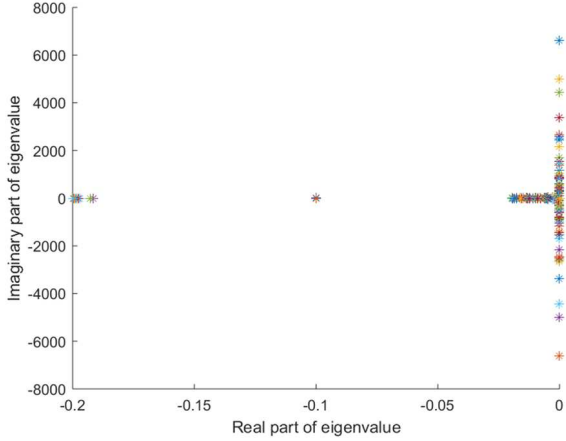

**Figure 1** Eigenvalues of UC with the 1st order model of turbine-governing system.

The eigenvalues in Fig. 1 are close to the real axis but with negative real parts. Fig. 2 shows the eigenvalues when both (1e) and (1f) are included. There are several real positive eigenvalues. This suggests that the closed-loop system (1), (6) under UC will be unstable. In order to improve stability of UC we propose modifications with the aim to approximately "decouple" the physical system (1) and the control system (6). We first derive the modified controller (10) and then argue that the approximate decoupling enhances stability.

In order to improve stability of UC we propose a modification with the aim to decouple the physical system (1) and control system (6). Observe that UC (6) is driven by state variables $\dot{\omega}_i, \omega_i, P_{ij}$ of the physical system through (6b), (6c), (6h). Other equations involve only control system quantities. To approximately decouple these physical quantities from the control system, we will use estimates of the disturbances $r_i$. They are unknown, however a method to approximate them is used. When the systems are decoupled, stability of the control system ensures stability of the physical system (1) regardless of the values of its parameters. If decoupled, the systems do not have feedback loop, therefore when control system reaches its equilibrium, it remains in it even when physical system has not yet converged to the desired point.

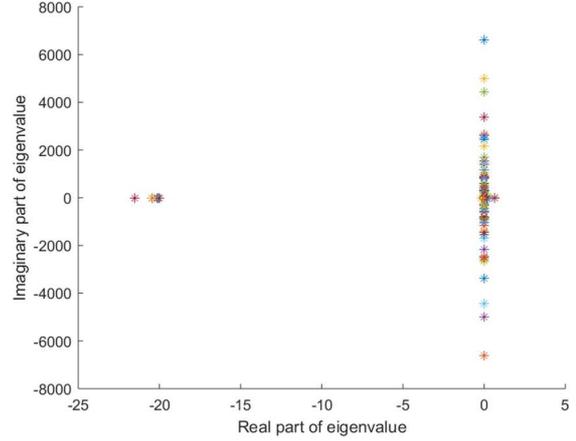

**Figure 2** Eigenvalues of UC with the 2nd order model of turbine-governing system. NB: the x-axis scale is different than in Fig. 1.

Let us consider load buses first. From (1c) it is possible to get estimation of disturbance:

$$r_i = D_i\omega_i + \sum_{ij\in\mathcal{E}} P_{ij} - p_i, i \in L \tag{11}$$

In case of generator buses (1b), we have the following estimations:

$$r_i + p_i^M = M_i\dot{\omega}_i + D_i\omega_i + \sum_{ij\in\mathcal{E}} P_{ij}, i \in L \tag{12}$$

Therefore equations (6b), (6c) are equivalent to

$$\dot{\lambda}_i = K_i^\lambda\left(-\sum_{ij\in\mathcal{E}} \hat{P}_{ij} - \alpha_i\lambda_i + r_i\right), i \in L, \tag{13a}$$

$$\dot{\lambda}_i = K_i^\lambda\left(-\sum_{ij\in\mathcal{E}} \hat{P}_{ij} + r_i + p_i^M\right), i \in G. \tag{13b}$$

The problem with (13b) is that the values of mechanical power injections $p_i^M$ are unknown. Therefore for generators, equations (10b), (10i), (10j) are used instead of equations (6b). Time constants $\tilde{T}^t$ and $\tilde{T}^g$ are estimates of true values $T^t$ and $T^g$ while $\tilde{p}_i^M$ and $\tilde{v}_i$ are estimates of $p_i^M$ and $v_i$. Equations (10i) and (10j) are used to approximate turbine and governor dynamics. Substituting (11) to (10c), and (12) to (11b) we get

$$\dot{\lambda}_i = K_i^\lambda\left(-\sum_{ij\in\mathcal{E}} \hat{P}_{ij} - \alpha_i\lambda_i + r_i\right), i \in L. \tag{14a}$$

$$\dot{\lambda}_i = K_i^\lambda\left(-\sum_{ij\in\mathcal{E}} \hat{P}_{ij} + r_i + p_i^M - \alpha_i\lambda_i - \tilde{p}_i^M\right), \tag{14b}$$

$$i \in G.$$

If the actual values $T^t$ and $T^g$ are approximated reasonably accurately by $\tilde{T}^t$ and $\tilde{T}^g$ respectively, then $p_i^M \approx \tilde{p}_i^M$ therefore (14b) is approximately equivalent to

$$\dot{\lambda}_i = K_i^\lambda \left( -\sum_{ij \in \varepsilon} \hat{P}_{ij} + r_i - \alpha_i \lambda_i \right), i \in G. \qquad (15)$$

Equations (15), similarly to (14a), do not depend on physical system variables directly. The control scheme uses $(\theta_i, \omega_i, P_{ij}, p_i^M, v_i)$ only to approximate the size of the disturbance $r_i$. As a result control scheme (10) is self-contained and it acts as an open loop system. Therefore the lag introduced by turbine and governor equations (1e), (1f) is removed from the system, which makes it stable.

Since in equilibrium $-\alpha_i \lambda_i^* - (\tilde{p}_i^M)^* = 0$, the modified control converges to the same desired point as original Unified Control (6). Additionally, if $\tilde{p}_i^M = p_i^M$ exactly, then decoupling allows us to exclude instabilities connected with introduction of governor dynamics equation (1f). Therefore stability proof derived for original UC is applicable for decoupled UC. As a result variables $\lambda_i$ are asymptotically stable and converge to the desired values $\lambda_i^*$. If variables $\tilde{p}_i^M$ approximate $p_i^M$ with some error $\tilde{p}_i^M - p_i^M = \sigma_i$, then $\lambda_i$ remain stable due to initial asymptotic stability.

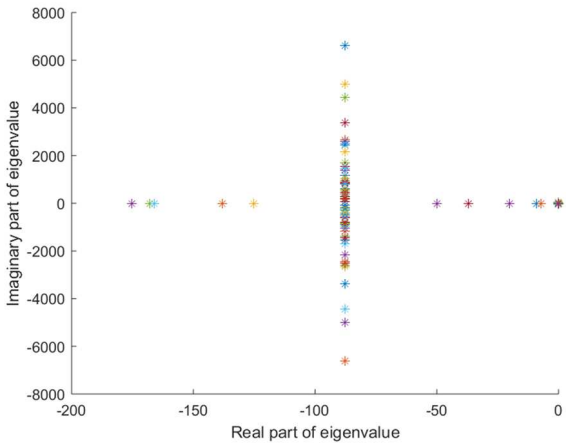

**Figure 3** Eigenvalues of the decoupled UC with the 2nd order model of turbine-governing system. NB: the x-axis scale is different than in Fig. 1 and Fig. 2.

Fig. 3 shows the eigenvalues for the decoupled UC control calculated for the IEEE 39-bus network with linearized power flow and with $\tilde{T}^t = T^t$ and $\tilde{T}^g = T^g$. Clearly all the eigenvalues that were close to the imaginary axis in Fig. 2 have been moved left and there are no eigenvalues with positive real parts.

## 6 Results of time-domain simulations

We have validated our modification by time-domain simulations using Power System Toolbox (PST) simulation program [3] with a full AC network model, a subtransient generator model, IEEE DC1 excitation model and IEEEGOV1 turbine governor model. System data was taken from [5]. For simplicity we did not implement area control (8e). Also, to compare UC/DUC with AGC on a fair basis, we assumed no load-side frequency control in UC/DUC – i.e. while $\lambda_i$ and $\phi_i$ were still calculated at load buses and their values exchanged with the neighboring buses, load demand itself was not changed. If load-side control was implemented, UC/DUC would obviously be even more effective.

Fig. 4 shows the plots of the control signal $\lambda_i$ at generator bus 34 for both the original and decoupled UC. We simulated a step change at node 38, $r_{38} = -7.35\ pu$. The control gains for UC were chosen to be $K_i^\lambda = 0.0477$, $K_i^\phi = 37.6991$ and $K_{ij}^{\rho^+} = K_{ij}^{\rho^-} = 0.0013$. The values of participation factors for UC were taken to be $\alpha_i = 1/R = 20$, where $R$ is the effective gain of the governing system so that the standard droop control was implemented.

The blue line in Fig. 4 shows the response of UC exhibiting undamped oscillations which, if not removed, would affect adversely the valves. The red line shows the response for the decoupled UC which is quite smooth.

Fig. 5 shows comparison of the frequency response for UC, Decoupled UC (DUC) and the standard AGC. The frequency response of UC and DUC is similar and much better than that of AGC in terms of the settling time as the control gains of UC and DUC have been chosen high. Trying to achieve a similar aggressive frequency control using AGC was impossible as it would have resulted in an unstable system. This shows another advantage of DUC – the control system is more stable than AGC and therefore the control can be quite aggressive. Note that frequency response of UC is smooth despite the oscillations in the control signal shown in Fig. 4, as the control signal acts on the physical system model which is effectively a low-pass filter.

Fig. 6 shows the values of flows on one of the congested lines. AGC is "blind" to the congestion so the green line stays above the line limit (showed by the black dashed line) following a disturbance. Decoupled UC behaves very similarly to the standard UC and drives the power flow to be within the line limits.

Note that in PST simulations the model of the turbine and its governing system is much more complicated than the assumed 2nd order model (1e), (1f) and its DUC emulator (10i), (10j). Hence the simulations stress-tested the assumption that the DUC emulates accurately the turbine and its governor in (10i), (10j), and that therefore $p_i^M \approx \tilde{p}_i^M$. To test further that assumption, Fig. 7 illustrates robustness of DUC to modelling errors in $\tilde{T}^t$ and $\tilde{T}^g$ when the turbine and its governor were modelled by a 2nd order system (1e), (1f). The blue line shows the results when the assumed and actual values of the parameters were the same while the red/green lines show the results when the assumed values of the parameters were 2 times bigger/smaller than the actual ones. Clearly such big parameter errors affect adversely the settling time or the overshoot but do not cause oscillations in the control signal seen in UC response in Fig.4.

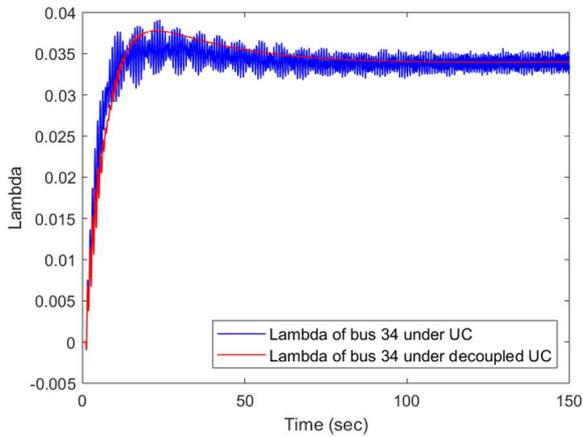

**Figure 4** Comparison of control signal $\lambda_i$ for generator at bus 34.

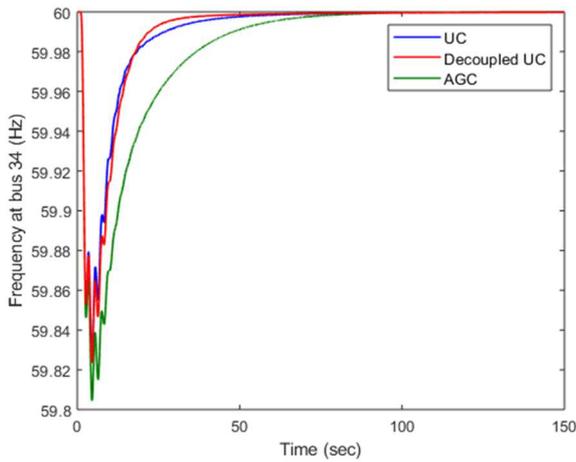

**Figure 5** Comparison of frequency plots for AGC, UC and decoupled UC.

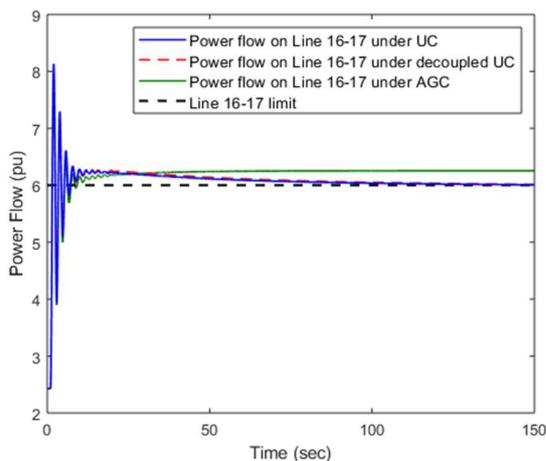

**Figure 6** Power flow on a congested line due to UC, modified UC and AGC.

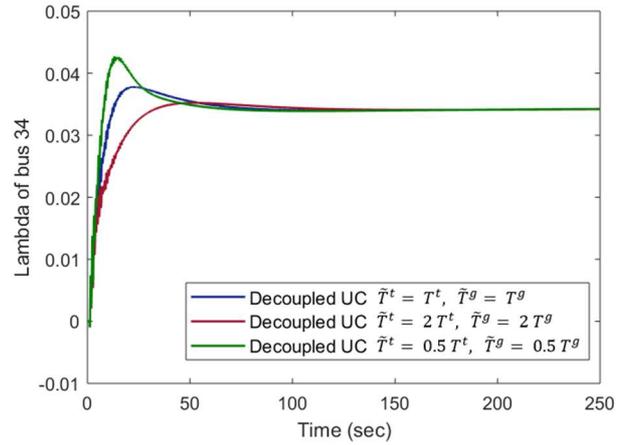

**Figure 7** Comparison of control signal for DUC with the presence of errors in $\tilde{T}^t$ and $\tilde{T}^g$.

# 7 Conclusions

This paper has analyzed stability of Unified Controller [1] that combines frequency control and congestion management therefore enabling corrective power system control. Earlier work by the authors of UC proved asymptotic stability of the methodology but the proof was based on a simplified first-order model of the turbine and turbine governor. We show that a higher order model of the turbine governor with high control gains leads to eigenvalues with small but positive real parts leading to sustained oscillations. Consequently we have developed a modification of the methodology that decouples the physical and control systems and therefore moves all the eigenvalues to the left and results in all the eigenvalues having negative real parts. Testing the modification on 39-bus model of New England power system containing 10 generators and using realistic dynamic models of the generators, exciters and turbine governors implemented in PST has proved the effectiveness of the methodology. The modification also allows for a more aggressive frequency control when compared to the standard AGC.